\documentclass{SLC}
\articlenumber{3}   
\graphicspath{{.}{../}{../Figures}{../../}}
\usepackage{tikz}

\newcommand*\circled[1]{\tikz[baseline=(char.base)]{
            \color{black}\draw (0,0.25) -- (0,0.75);
            \draw (1,0.25) -- (1,0.5);
            \draw (2,0.5) -- (2,0.75);
            \draw (-0.25,0.5) -- (1,0.5);
            \draw (2,0.5) -- (2.25,0.5);
            \draw [black, fill=black] (1,1) circle(2pt);
            \draw [black, fill=black] (2,0) circle(2pt);
            \node[shape=circle,draw,inner sep=2pt] at (0,0) (char) {#1};
            \node[shape=circle,draw,inner sep=2pt] at (0,1) (char) {#1};
            \node[shape=circle,draw,inner sep=2pt] at (1,0) (char) {#1};
            \node[shape=circle,draw,inner sep=2pt] at (2,1) (char) {#1};
            }}
\newcommand*\circledtwo[2]{\tikz[baseline=(char.base)]{
            \draw (2,0.25) -- (2,0.5);
            \draw (0,0.5) -- (0,0.75);
            \draw (0,0.5) -- (2,0.5);
            \draw [black, fill=black] (0,0) circle(2pt);
            \draw [black, fill=black] (1,1) circle(2pt);
            \draw [black, fill=black] (2,1) circle(2pt);
            \node[shape=circle,draw,inner sep=2pt] at (1,0) (char) {#1};
            \node[shape=circle,draw,inner sep=2pt] at (2,0) (char) {#2};
            \node[shape=circle,draw,inner sep=2pt] at (0,1) (char) {#2};
            }}
\newcommand*\circledthree[2]{\tikz[baseline=(char.base)]{
            \draw (2,0.25) -- (2,0.75);
            \draw [black, fill=black] (0,0) circle(2pt);
            \draw [black, fill=black] (1,1) circle(2pt);
            \draw [black, fill=black] (0,1) circle(2pt);
            \node[shape=circle,draw,inner sep=2pt] at (1,0) (char) {#1};
            \node[shape=circle,draw,inner sep=2pt] at (2,0) (char) {#2};
            \node[shape=circle,draw,inner sep=2pt] at (2,1) (char) {#2};
            }}
\newcommand*\circledfour[2]{\tikz[baseline=(char.base)]{
            \draw (0,0.5) -- (0,0.25);
            \draw (2,0.5) -- (2,0.75);
            \draw (-0.25,0.5) -- (0,0.5);
            \draw (2.25,0.5) -- (2,0.5);
            \draw [black, fill=black] (2,0) circle(2pt);
            \draw [black, fill=black] (1,1) circle(2pt);
            \draw [black, fill=black] (0,1) circle(2pt);
            \node[shape=circle,draw,inner sep=2pt] at (1,0) (char) {#1};
            \node[shape=circle,draw,inner sep=2pt] at (0,0) (char) {#2};
            \node[shape=circle,draw,inner sep=2pt] at (2,1) (char) {#2};
            }}
\newcommand*\circledfive[2]{\tikz[baseline=(char.base)]{
            \draw (0,0.75) -- (0,0.25);
            \draw [black, fill=black] (2,0) circle(2pt);
            \draw [black, fill=black] (1,1) circle(2pt);
            \draw [black, fill=black] (2,1) circle(2pt);
            \node[shape=circle,draw,inner sep=2pt] at (1,0) (char) {#1};
            \node[shape=circle,draw,inner sep=2pt] at (0,0) (char) {#2};
            \node[shape=circle,draw,inner sep=2pt] at (0,1) (char) {#2};
            }}
\newcommand*\thedasep[1]{\tikz[baseline=(char.base)]{
            \draw (0.5,10.2) node[minimum size=8.0cm,draw,circle] {};
            \draw (0.5,17.7) node[minimum size=5.21cm,draw,circle] {};
            \draw (0.5,2.7) node[minimum size=5.21cm,draw,circle] {};

            \begin{scope}[ color                 = black
                 , execute at begin node = $\displaystyle
                 , execute at end node   = $]
                  \node at (-4.8,10.2) {{\textnormal\ASEP}(2,1,0)} ;
                  \node at (-3.4,17.7) {{\textnormal\ASEP}(1,1,0)} ;
                  \node at (-3.4,2.7) {{\textnormal\ASEP}(2,2,0)} ;
            \end{scope}

            \begin{scope}[ color                 = blue]
            \draw [<->] (-1.7305758785247998,10.171709811667625) -- (-1.7305758785247998,10.571709811667626);
            \draw [<->] (-0.590787917055414,12.145880469772806) -- (-0.937198079055414,11.945880469772806);
            \draw [<->] (1.5907879170554138,12.145880469772806) -- (1.9371980790554137,11.945880469772806);
            \draw [<->] (2.7305758785248,10.171709811667625) -- (2.7305758785248,10.571709811667626);
            \draw [<->] (1.6397879615661848,8.28240971882309) -- (1.2933777995661848,8.082409718823088);
            \draw [<->] (-0.639787961566185,8.28240971882309) -- (-0.2933777995661849,8.082409718823088);
            \draw [<->] (1.6847227529,9.515999999905187) -- (-0.3937382170999999,10.715999999905186);
            \draw [<->] (-0.6847227529000002,9.515999999905187) -- (1.3937382171000001,10.715999999905186);
            \draw [<->] (0.5,9.168000000405186) -- (0.5,11.568000000405185);
            \draw [<->] (0.08051323580980352,16.99425215426003) -- (0.42692339780980326,16.79425215426003);
            \draw [<->] (0.0985478192978495,18.416160117102656) -- (-0.2478623427021498,18.216160117102657);
            \draw [<->] (1.320938945374892,17.689587728919978) -- (1.320938945374892,18.08958772891998);
            \draw [<->] (0.9194867641901967,1.9942521542600304) -- (0.5730766021901967,1.7942521542600307);
            \draw [<->] (0.9014521807021505,3.416160117102657) -- (1.24786234270215,3.216160117102657);
            \draw [<->] (-0.3209389453748921,2.68958772891998) -- (-0.3209389453748921,3.08958772891998);
            \end{scope}

            \begin{scope}[ color                 = red]
            \draw [<->] (2.5907879170554138,12.515880469772806) -- (1.17692339780980326,15.79425215426003);
            \draw [<->] (2.5907879170554138,8.134119530227192) -- (1.97692339780980326,3.805747845739968);
            \end{scope}

            \draw (2.1547005383792506,2.411324865405187) -- (1.6547005383792506,2.411324865405187);
            \draw (2.2797005383792506,2.627831216351297) -- (2.0297005383792506,3.060843918243516);
            \draw (1.5297005383792506,2.627831216351297) -- (1.7797005383792506,3.060843918243516);
            \node[shape=circle,draw,inner sep=2pt] at (1.4047005383792506,2.411324865405187) (char) {0};
            \node[shape=circle,draw,inner sep=2pt] at (2.4047005383792506,2.411324865405187) (char) {2};
            \node[shape=circle,draw,inner sep=2pt] at (1.9047005383792506,3.277350269189626) (char) {2};

            \draw (0.047649730810373825,1.1948185144590775) -- (-0.4523502691896262,1.1948185144590775);
            \draw (0.17264973081037382,1.411324865405187) -- (-0.07735026918962618,1.8443375672974065);
            \draw (-0.5773502691896262,1.411324865405187) -- (-0.3273502691896262,1.8443375672974065);
            \node[shape=circle,draw,inner sep=2pt] at (-0.7023502691896262,1.1948185144590775) (char) {2};
            \node[shape=circle,draw,inner sep=2pt] at (0.2976497308103738,1.1948185144590775) (char) {0};
            \node[shape=circle,draw,inner sep=2pt] at (-0.20235026918962618,2.060843918243516) (char) {2};

            \draw (0.047649730810373825,3.6278312163512965) -- (-0.4523502691896262,3.6278312163512965);
            \draw (0.17264973081037382,3.8443375672974067) -- (-0.07735026918962618,4.2773502691896255);
            \draw (-0.5773502691896262,3.8443375672974067) -- (-0.3273502691896262,4.2773502691896255);
            \node[shape=circle,draw,inner sep=2pt] at (-0.7023502691896262,3.6278312163512965) (char) {2};
            \node[shape=circle,draw,inner sep=2pt] at (0.2976497308103738,3.6278312163512965) (char) {2};
            \node[shape=circle,draw,inner sep=2pt] at (-0.20235026918962618,4.493856620135736) (char) {0};

            \draw (-1.1547005383792506,17.411324865405184) -- (-0.6547005383792506,17.411324865405184);
            \draw (-1.2797005383792506,17.627831216351296) -- (-1.0297005383792506,18.060843918243513);
            \draw (-0.5297005383792506,17.627831216351296) -- (-0.7797005383792506,18.060843918243513);
            \node[shape=circle,draw,inner sep=2pt] at (-0.4047005383792506,17.411324865405184) (char) {0};
            \node[shape=circle,draw,inner sep=2pt] at (-1.4047005383792506,17.411324865405184) (char) {1};
            \node[shape=circle,draw,inner sep=2pt] at (-0.9047005383792506,18.277350269189625) (char) {1};

            \draw (0.9523502691896262,16.194818514459076) -- (1.4523502691896262,16.194818514459076);
            \draw (0.8273502691896262,16.411324865405184) -- (1.0773502691896262,16.844337567297405);
            \draw (1.5773502691896262,16.411324865405184) -- (1.3273502691896262,16.844337567297405);
            \node[shape=circle,draw,inner sep=2pt] at (1.7023502691896262,16.194818514459076) (char) {1};
            \node[shape=circle,draw,inner sep=2pt] at (0.7023502691896262,16.194818514459076) (char) {0};
            \node[shape=circle,draw,inner sep=2pt] at (1.2023502691896262,17.060843918243513) (char) {1};

            \draw (0.9523502691896262,18.627831216351296) -- (1.4523502691896262,18.627831216351296);
            \draw (0.8273502691896262,18.844337567297405) -- (1.0773502691896262,19.277350269189625);
            \draw (1.5773502691896262,18.844337567297405) -- (1.3273502691896262,19.277350269189625);
            \node[shape=circle,draw,inner sep=2pt] at (1.7023502691896262,18.627831216351296) (char) {1};
            \node[shape=circle,draw,inner sep=2pt] at (0.7023502691896262,18.627831216351296) (char) {1};
            \node[shape=circle,draw,inner sep=2pt] at (1.2023502691896262,19.493856620135734) (char) {0};

            \draw (0.25,7.411324865405186) -- (0.75,7.411324865405186);
            \draw (0.125,7.627831216351296) -- (0.375,8.060843918243515);
            \draw (0.875,7.627831216351296) -- (0.625,8.060843918243515);
            \node[shape=circle,draw,inner sep=2pt] at (1.0,7.411324865405186) (char) {0};
            \node[shape=circle,draw,inner sep=2pt] at (0.0,7.411324865405186) (char) {1};
            \node[shape=circle,draw,inner sep=2pt] at (0.5,8.277350269189625) (char) {2};

            \draw (-1.9150635094610964,8.661324865405186) -- (-1.4150635094610964,8.661324865405186);
            \draw (-2.0400635094610964,8.877831216351296) -- (-1.7900635094610964,9.310843918243515);
            \draw (-1.2900635094610966,8.877831216351296) -- (-1.5400635094610966,9.310843918243515);
            \node[shape=circle,draw,inner sep=2pt] at (-1.1650635094610966,8.661324865405186) (char) {0};
            \node[shape=circle,draw,inner sep=2pt] at (-2.1650635094610964,8.661324865405186) (char) {2};
            \node[shape=circle,draw,inner sep=2pt] at (-1.6650635094610966,9.527350269189625) (char) {1};

            \draw (2.4150635094610893,8.661324865405186) -- (2.9150635094610973,8.661324865405186);
            \draw (2.2900635094610973,8.877831216351296) -- (2.5400635094610973,9.310843918243515);
            \draw (3.0400635094610964,8.877831216351296) -- (2.7900635094610973,9.310843918243515);
            \node[shape=circle,draw,inner sep=2pt] at (3.1650635094610964,8.661324865405186) (char) {1};
            \node[shape=circle,draw,inner sep=2pt] at (2.1650635094610964,8.661324865405186) (char) {0};
            \node[shape=circle,draw,inner sep=2pt] at (2.6650635094610964,9.527350269189625) (char) {2};

            \draw (2.4150635094610893,11.161324865405186) -- (2.9150635094610973,11.161324865405186);
            \draw (2.2900635094610973,11.377831216351296) -- (2.5400635094610973,11.810843918243515);
            \draw (3.0400635094610964,11.377831216351296) -- (2.7900635094610973,11.810843918243515);
            \node[shape=circle,draw,inner sep=2pt] at (3.1650635094610964,11.161324865405186) (char) {2};
            \node[shape=circle,draw,inner sep=2pt] at (2.1650635094610964,11.161324865405186) (char) {0};
            \node[shape=circle,draw,inner sep=2pt] at (2.6650635094610964,12.027350269189625) (char) {1};

            \draw (-1.9150635094610964,11.161324865405186) -- (-1.4150635094610964,11.161324865405186);
            \draw (-2.0400635094610964,11.377831216351296) -- (-1.7900635094610964,11.810843918243515);
            \draw (-1.2900635094610966,11.377831216351296) -- (-1.5400635094610966,11.810843918243515);
            \node[shape=circle,draw,inner sep=2pt] at (-1.1650635094610966,11.161324865405186) (char) {1};
            \node[shape=circle,draw,inner sep=2pt] at (-2.1650635094610964,11.161324865405186) (char) {2};
            \node[shape=circle,draw,inner sep=2pt] at (-1.6650635094610966,12.027350269189625) (char) {0};

            \draw (0.25,12.411324865405186) -- (0.75,12.411324865405186);
            \draw (0.125,12.627831216351296) -- (0.375,13.060843918243515);
            \draw (0.875,12.627831216351296) -- (0.625,13.060843918243515);
            \node[shape=circle,draw,inner sep=2pt] at (1.0,12.411324865405186) (char) {2};
            \node[shape=circle,draw,inner sep=2pt] at (0.0,12.411324865405186) (char) {1};
            \node[shape=circle,draw,inner sep=2pt] at (0.5,13.277350269189625) (char) {0};

            }}

\author[D.~Ash]{David W. Ash\textsuperscript{1}\protect\orcid{0000-0001-8970-5784}}
\address{\textsuperscript{1}OptIn Inc., San Francisco Bay Area, CA, USA; \website{https://orcid.org/0000-0001-8970-5784}}

\title[Introducing DASEP: the doubly asymmetric simple exclusion process]{Introducing DASEP:\\ the doubly asymmetric simple exclusion process }

\def\ASEP{{\footnotesize ASEP}}
\def\DASEP{{\footnotesize DASEP}}

\abstract{Research in combinatorics has often explored the asymmetric simple exclusion process ({\ASEP}). The {\ASEP}, inspired by
examples from statistical mechanics, involves particles of various species moving around a lattice. With the 
traditional {\ASEP} particles of a given species can move but do not change species. In this paper a new combinatorial
formalism, the {\DASEP} (doubly asymmetric simple exclusion process), is explored. The {\DASEP} is inspired by biological
processes where, unlike the {\ASEP}, the particles can change from one species to another. The combinatorics of the 
{\DASEP} on a one dimensional lattice are explored, including the associated generating function. The stationary 
probabilities of the {\DASEP} are explored, and results are proven relating these stationary probabilities to those of 
the simpler {\ASEP}.}

 
\begin{document}

\maketitle

\section{Introduction}

The {\ASEP} (asymmetric simple exclusion process) is a structure that has frequently been referred to in the combinatorics
literature. In its simplest form, the {\ASEP} consists of a one dimensional infinite lattice, with each point on the
lattice being populated with either a particle or a hole. At random intervals, each particle attempts to move either to
the left or the right with different but fixed probabilites (hence the term `asymmetric'). The {\ASEP} can be thought of
as a form of Markov process as noted in~\cite{multiline} by Corteel et al.~ Multiline queues~\cite{ferrari} were
introduced by Ferrari et al.~ as a combinatorial approach to the analysis of the {\ASEP}. Originally the {\ASEP} particles
were thought of as all belonging to a single species. More recent work by Cantini et al.~\cite{cantini} generalized
the concept to multiple species and uncovered a link with Macdonald polynomials. Although we focus on the
homogeneous {\ASEP} (transition probabilities do not depend on position in the lattice), several 
researchers (Lam et al.~\cite{lam}, Ayyer et al.~\cite{ayyer}, Cantini~\cite{cantiniarc}, Mandelshtam~\cite{mandelshtam},
and Kim et al.~\cite{kim}) have
explored the inhomogeneous {\ASEP} in which transition probabilities do depend on lattice position.

\section{Definitions}

Following~\cite{multiline}, a partition can be defined as follows:

\begin{definition} \label{defn:name2}
A partition $ \lambda $ is a nonincreasing sequence of $ n $ nonnegative integers
$ \lambda = (\lambda_1 \ge \lambda_2 \ge ... \ge \lambda_n \ge 0)$.
\end{definition}

We will start by working through a simple example of
the {\ASEP} before introducing the new concept of the {\DASEP}. We will ordinarily write a partition as defined above as an $n$-tuple:
$\lambda=(\lambda_1,\lambda_2,...,\lambda_n)$. 

\begin{definition}\label{SnL}
We write $S_n(\lambda )$ to mean the set of all permutations of $\lambda$. 
\end{definition}

\begin{example}
So,
for $\lambda=(2,2,1)$, $S_3(\lambda)=\{(2,2,1),(2,1,2),(1,2,2)\}$. 
\end{example}

The multispecies asymmetric simple exclusion
process {\ASEP}($\lambda$) is then defined to be a Markov process on $S_n(\lambda )$ with certain specific transition probabilities:

\begin{definition}
For all partitions $\lambda$ as defined in Definition \ref{defn:name2},
ASEP($\lambda$) is a Markov process on $S_n(\lambda )$. We let $t$ be a constant with $0\le t\le 1$. 
The transition probability $P_{\mu,\nu}$ between 
two permutations $\mu\in S_n(\lambda )$ and $\nu\in S_n(\lambda )$ is given by:
\begin{itemize}
  \item If $\mu=(\mu_1,\dots,\mu_k,i,j,\mu_{k+2},\dots,\mu_n)$ and 
$\nu=(\mu_1,\dots,\mu_k,j,i,\mu_{k+2},\dots,\mu_n)$, with $i\ne j$,  then $P_{\mu,\nu}=\frac{t}{n}$ if
$i>j$ and $P_{\mu,\nu}=\frac{1}{n}$ if $j>i$.
  \item If $\mu=(i,\mu_2,\mu_3,\dots,\mu_{n-1},j)$ and
$\nu=(j,\mu_2,\mu_3,\dots,\mu_{n-1},i)$ with $i\ne j$, then $P_{\mu,\nu}=\frac{t}{n}$ if
$j>i$ and $P_{\mu,\nu}=\frac{1}{n}$ if $i>j$.
  \item If neither of the above conditions apply but $\nu\ne\mu$ then $P_{\mu,\nu}=0$. If $\nu=\mu$ then 
$P_{\mu,\mu}=1-\sum_{\nu\ne\mu}P_{\mu,\nu}$.
\end{itemize}
\end{definition}

It is possible to compute steady state probabilities for {\ASEP}($\lambda$). For the purposes of the example that we will
develop as we introduce {\DASEP}, we are primarily interested in {\ASEP}($\lambda$) for $\lambda=(2,2,0)$, $\lambda=(2,1,0)$, 
and $\lambda=(1,1,0)$, so we will focus mostly on these three processes as we work through the computation of the
steady state probabilities. Continuing to follow~\cite{multiline} as we develop this example, to compute these 
probabilities we need to define the concept of a multiline queue. 

\begin{definition}
A ball system $B$ is an $L\times n$
matrix each element of which is either $0$ or~$1$. Moreover for all $i$ the number of $1$'s in row $i+1$ is less than
or equal to the number of $1$'s in row $i$.
\end{definition}

\begin{definition}
Given a ball system $B$ a multiline queue $Q$ is obtained by augmenting $B$ with a labeling and matching system.
Each cell in $B$ will be labelled with a number from $0$ to $L$ inclusive, and each cell with a $1$ element in row
$i+1$, for $i\ge 1$, will be matched to a cell with a $1$ element in row $i$. Such a matching must be obtained through
an application of the following algorithm:

\begin{itemize}
  \item Step 1: Find the highest numbered row with unlabelled $1$ elements. Label each of those elements
with the number of the row. If this is row $1$, or there are no remaining unlabelled $1$ elements in the matrix,
exit.
  \item Step 2: Find the row with labelled but unmatched elements. If this is row $1$, go back to step 1.
If it is row $i+1$, for $i\ge 1$, first match each labelled but unmatched element that can be matched to an unlabelled element
directly below it to that element. This is considered a \textit{trivial} match. Then proceed from right to left
(highest to lowest numbered columns) matching each remaining labelled but unmatched element to an unlabelled element
in the row below--these are the \textit{nontrivial} matches. Give all newly matched elements in row $i$ the same label as the element
it has just been matched to. Repeat step 2.
\end{itemize}
\end{definition}

A multiline queue is often visualized as a ball system with an element with a $1$ value being shown as a ball and a $0$
value by the absence of a ball. Matches between elements (balls) are drawn by lines between the matched balls.
The following shows a multiline queue associated with {\ASEP}($\lambda$) where $\lambda=(2,2,0)$. Note that the line
matching the ball at upper right to the one at the lower middle wraps around to the right.
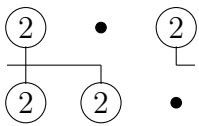
\begin{center}
\circled{2}
\end{center}

The labels in the bottom row determine the partition of the associated {\ASEP}. The above multiline queue has $\lambda=(2,2,0)$
since the bottom row includes two $2$'s and a $0$--by convention an element without a ball is assumed to be labeled with
a $0$. Likewise the following would be a multiline queue with $\lambda=(2,1,0)$:
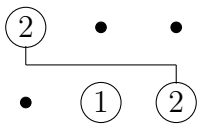
\begin{center}
\circledtwo{1}{2}
\end{center}

Each multiline queue is also associated with a permutation $\alpha\in S_n(\lambda)$ corresponding to the labels of its
bottom row in \textit{unsorted} order. For example, for the above multiline queue, $\lambda=(2,1,0)$ but
$\alpha=(0,1,2)$. We will write $\lambda(Q)=\lambda$ and $\alpha(Q)=\alpha$.

\section{Steady state probabilities with example}

To determine steady state probabilities--and continue with the example started in the introduction--we next assign to 
each nontrivial matching $p$ in $Q$ two values $f(p)$ and 
$s(p)$. $f(p)$ is the number of choices that were available for the match when the match was made. $s(p)$ is the number
of legal matches that were skipped, if we imagine ourselves considering possible matches from left to right and wrapping
around the end if needed, before the actual choice was made. We can then define a weight on $p$ as
$\textnormal{wt}(p)=\frac{(1-t)t^{s(p)}}{1-t^{f(p)}}$. Here we are proceeding from~\cite{multiline} but with the
simplifying assumption that $q=1$, since in the sequel we will rely on a theorem that requires $q=1$. Next we can define
a weight on the entire multiline queue $\textnormal{wt}(Q)=\prod_{p\in Q} \textnormal{wt}(p)$ where the product is taken
over all \textit{nontrivial} matches $p$ in $Q$. A theorem due to Martin~\cite{martin} then gives the required steady
state probabilities:
\begin{equation}
 \textnormal{Pr}(\alpha)=\frac{\sum_{\alpha(Q)=\alpha} \textnormal{wt}(Q)}{\sum_{\lambda(Q)=\lambda} \textnormal{wt}(Q)}. \end{equation}

Before moving on to the {\DASEP}, we need to evaluate the steady state probabilities for the examples that we will ultimately
use to develop the {\DASEP}. For the above multiline queue, there is exactly one nontrivial pair $p$. When this pair is
matched, there are two available options so $f(p)=2$. As we picked the second available option, $s(p)=1$. 
So $\textnormal{wt}(Q)=\frac{(1-t)t}{1-t^2}$. As noted above, $\alpha=(0,1,2)$ and the only other multiline queue with
$\alpha=(0,1,2)$ is as follows:

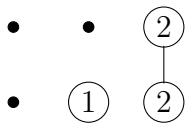
\begin{center}
\circledthree{1}{2}
\end{center}

Here there is no nontrivial matching pair, so $\textnormal{wt}(Q)=1$. Hence:
\begin{equation} \sum_{\alpha(Q)=(0,1,2)} \textnormal{wt}(Q) = 1 + \frac{(1-t)t}{1-t^2} = \frac{1+2t}{1+t}. \end{equation}

For reasons of symmetry:
\begin{equation} \sum_{\alpha(Q)=(0,1,2)} \textnormal{wt}(Q) = \sum_{\alpha(Q)=(1,2,0)} \textnormal{wt}(Q) = 
\sum_{\alpha(Q)=(2,0,1)} \textnormal{wt}(Q) = \frac{1+2t}{1+t}. \end{equation}

Next we look at $\alpha=(0,1,2)$, for which there are also two multiline queues. The first of these is as follows:
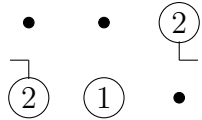
\begin{center}
\circledfour{1}{2}
\end{center}

Here there is one nontrivial matching pair $p$. When this pair is
matched, there are two available options so $f(p)=2$. As we picked the first available option, $s(p)=0$. So 
$\textnormal{wt}(Q)=\frac{1-t}{1-t^2}$. The other multiline queue with
$\alpha=(2,1,0)$ is as follows:

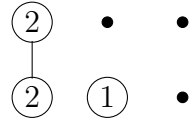
\begin{center}
\circledfive{1}{2}
\end{center}

Again there is no nontrivial matching pair, so $\textnormal{wt}(Q)=1$. Hence:
\begin{equation} \sum_{\alpha(Q)=(2,1,0)} \textnormal{wt}(Q) = 1 + \frac{1-t}{1-t^2} = \frac{2+t}{1+t}.\end{equation}
For reasons of symmetry, one has
\begin{equation} \sum_{\alpha(Q)=(2,1,0)} \textnormal{wt}(Q) = \sum_{\alpha(Q)=(1,0,2)} \textnormal{wt}(Q) = 
\sum_{\alpha(Q)=(0,2,1)} \textnormal{wt}(Q) = \frac{2+t}{1+t}. \end{equation}
So we get
\begin{equation} \sum_{\lambda(Q)=(2,1,0)} \textnormal{wt}(Q) = 3(\frac{1+2t}{1+t})+3(\frac{2+t}{1+t}) = 9.\end{equation}
We are now ready to give the steady state probabilities
\begin{equation} \textnormal{Pr}(0,1,2) = \textnormal{Pr}(1,2,0) = \textnormal{Pr}(2,0,1) = \frac{1+2t}{9(1+t)} \end{equation}
and
\begin{equation} \textnormal{Pr}(2,1,0) = \textnormal{Pr}(1,0,2) = \textnormal{Pr}(0,2,1) = \frac{2+t}{9(1+t)}. \end{equation}
Trivial computations also give
\begin{equation} \textnormal{Pr}(0,1,1) = \textnormal{Pr}(1,1,0) = \textnormal{Pr}(1,0,1) = \frac{1}{3}\end{equation}
and
\begin{equation} \textnormal{Pr}(0,2,2) = \textnormal{Pr}(2,2,0) = \textnormal{Pr}(2,0,2) = \frac{1}{3}. \end{equation}

This concludes our computation for the steady state probabilities of this model;
in the next section we introduce the {\DASEP} model.
\pagebreak
\begin{figure}
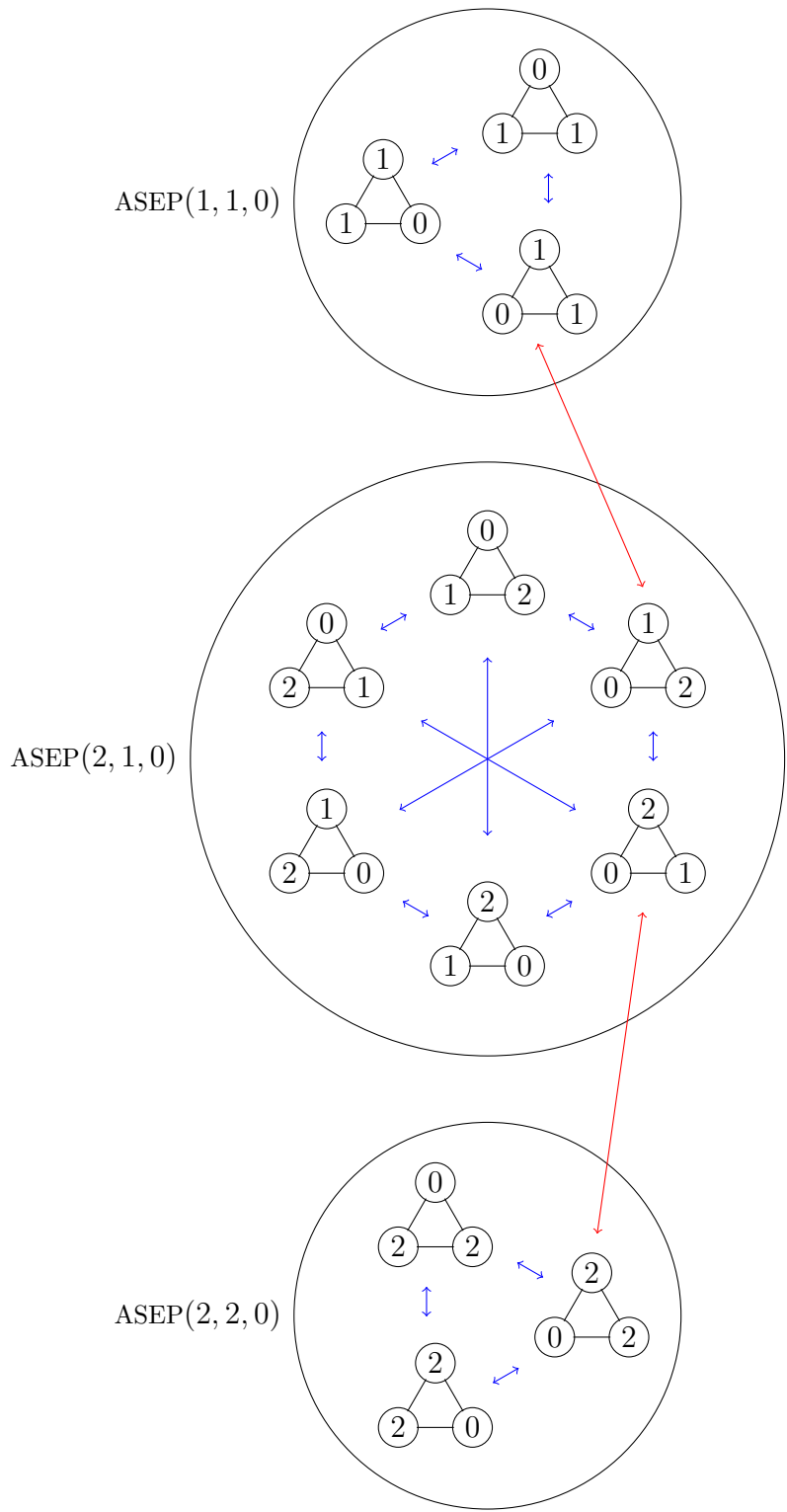
\thedasep{0}
\setlength\abovecaptionskip{2mm}
\setlength\belowcaptionskip{-4mm}
\caption{An example of {\DASEP}($n,p,q$): {\DASEP}(3,2,2).}\label{fig1} 
\end{figure}

\qquad\vspace{-9mm}

Each of the $12$ small triangles represents a state 
of the {\DASEP} in the circular lattice with $n=3$ sites with $q=2$ balls (i.e., the nonzero labels),
each nonzero label is $\leq p=2$ (i.e., one has $p=2$ species). Each state corresponds to a permutation of the partition $(2,2,0)$, $(2,1,0)$, or $(1,1,0)$.  
The transitions are explained in Definition~\ref{defn:name1} hereafter.
\pagebreak

\section{Doubly asymmetric simple exclusion process}
\label{section:dasep-start}

We are now ready to introduce the {\DASEP} (doubly asymmetric simple exclusion process). 
While the {\ASEP} is
inspired by statistical mechanics where particles do not change species, the {\DASEP}, by contrast, is inspired by
biological processes where particles~can change species,
which we denote 
by {\DASEP}$(n,p,q)$ where $n$ is the number
of positions on the lattice, $p$ is the number of types of species, and $q$ is the number of particles. 

\begin{definition} \label{defn:name1}
For all positive integers $n$, $p$, and $q$ with $n>q$, {\DASEP}$(n,p,q)$
is a Markov process on the set
$\displaystyle{  \bigcup_{\substack{\lambda_1 \le p, \ \lambda^\prime_1 = q}} S_n(\lambda)}$,
where one uses the notation of Definition~\ref{SnL}, and where $\lambda^\prime_1 = q$ refers to the dual partition~\cite{macdonald} of $\lambda$, namely $\lambda^\prime$,
and uses the fact that $\lambda^\prime_1$ gives the number of nonzero terms in the original partition $\lambda$.
The transition probability $P_{\mu,\nu}$ on two permutations $\mu$ and $\nu$
is as follows:

\begin{itemize}
  \item If $\mu=(\mu_1,\dots,\mu_k,i,j,\mu_{k+2},\dots,\mu_n)$ and 
$\nu=(\mu_1,\mu_2,...,\mu_k,j,i,\mu_{k+2},\dots,\mu_n)$ with $i\ne j$, then $P_{\mu,\nu}=\frac{t}{3n}$ if
$i>j$ and $P_{\mu,\nu}=\frac{1}{3n}$ if $j>i$.
  \item If $\mu=(i,\mu_2,\dots,\mu_{n-1},j)$ and
$\nu=(j,\mu_2,\mu_3,...,\mu_{n-1},i)$ with $i\ne j$, then $P_{\mu,\nu}=\frac{t}{3n}$ if
$j>i$ and $P_{\mu,\nu}=\frac{1}{3n}$ if $i>j$.
  \item If $\mu=(\mu_1,\dots,\mu_k,i,\mu_{k+2},\dots,\mu_n)$ and 
$\nu=(\mu_1,\dots,\mu_k,i+1,\mu_{k+2},\dots,\mu_n)$ with $i\ge 1$, then $P_{\mu,\nu}=\frac{u}{3n}$.
  \item If $\mu=(\mu_1,\dots,\mu_k,i+1,\mu_{k+2},\dots,\mu_n)$ and 
$\nu=(\mu_1,\dots,\mu_k,i,\mu_{k+2},\dots,\mu_n)$ with $i\ge 1$, then $P_{\mu,\nu}=\frac{1}{3n}$.
  \item If none of the above conditions apply but $\nu\ne\mu$ then $P_{\mu,\nu}=0$. If $\nu=\mu$ then 
$P_{\mu,\mu}=1-\sum_{\nu\ne\mu}P_{\mu,\nu}$.
\end{itemize}
\end{definition}
\smallskip

Figure~\ref{fig1} shows the simple example of the {\DASEP} that we are working through.
All possible transitions within a single {\ASEP} 
(the first and second bullet points in the definition above) are shown with blue arrows on this diagram.  
To keep the diagram relatively clean in appearance, only selected transitions between different
ASEPs (the third and fourth bullet points) are shown (with red arrows). Other {\ASEP}s such as {\ASEP}$(1,0,0)$
or {\ASEP}$(2,0,0)$ are not shown since these are not part of {\DASEP}$(3,2,2)$. This is because, per
Definition~\ref{defn:name1}, for {\DASEP}$(3,2,2)$ we always have $\lambda^\prime_1 = 2$, whereas for
ASEP$(1,0,0)$ and {\ASEP}$(2,0,0)$, we would have $\lambda^\prime_1 = 1$.
\smallskip

Similar to with the {\ASEP}, with the {\DASEP} we wish to compute steady state probabilities for permutations $\alpha$
which we will call $\textnormal{Pd}(\alpha)$. We will focus on continuing to develop the example we have been working
on which turns out to be {\DASEP}$(3,2,2)$. Here $n=3$ means that the particles move on the the circular lattice with $3$ sites,
$p=2$ means that each particle is allowed to take on the value $0$, $1$, or $2$, and $q=2$ means that each 
permutation~$\alpha$ has exactly $2$ nonzero values. We therefore find ourselves interested in the following $12$
steady state probabilities:  
\begin{align}&\textnormal{Pd}(0,1,1), \textnormal{Pd}(0,1,2), \textnormal{Pd}(0,2,1), \textnormal{Pd}(0,2,2),
   \textnormal{Pd}(1,0,1), \textnormal{Pd}(1,0,2), \\
&\textnormal{Pd}(1,1,0), \textnormal{Pd}(1,2,0), \textnormal{Pd}(2,0,1), \textnormal{Pd}(2,0,2),
   \textnormal{Pd}(2,1,0), \textnormal{Pd}(2,2,0). \end{align}
\pagebreak

Note here that particles in the {\DASEP} are allowed to switch back and forth between species $1$ and $2$, but not
back and forth from $0$ to anything else. That is because a value of $0$ is understood to not so much be a species
but the absence of a species. Due to symmetries we can now focus on solving for the following four probabilities:
\begin{equation} w=\textnormal{Pd}(0,1,1), \quad x=\textnormal{Pd}(0,1,2), \quad  y=\textnormal{Pd}(0,2,1), \quad  z=\textnormal{Pd}(0,2,2). \end{equation}

From the above transition probabilities, this reduces to solving the system
\begin{equation}\begin{cases}
 2uw = x+y \\
 (2+t)x + x + ux = (1+2t)y+z+uw \\
 (1+2t)y + y + uy = (2+t)x + uw + z \\
2z = u(x+y)
\end{cases} \end{equation}
which in turn implies the relation
\begin{equation} (5+2t+u)x = (3+4t+u)y. \end{equation}

We can then ask ourselves the question of when the proportions of steady state probabilities for the {\DASEP}
are the same as for the previous {\ASEP}. Noting that
$\textnormal{Pr}(0,1,2) = \frac{1+2t}{9(1+t)}$ and
$\textnormal{Pr}(2,1,0) = \frac{2+t}{9(1+t)}$ such equality will happen if
$(5+2t+u)(1+2t) = (3+4t+u)(2+t)$,
or
$5+2t+u+10t+4t^2+2tu = 6+8t+2u+3t+4t^2+tu$,
or $t(1+u) = 1+u.$
So this will happen iff $t=1$. We have therefore proven the following proposition.
\begin{proposition}
If $D=\textnormal{{\DASEP}}(3,2,2)$ is parameterized as described above by $t$ and $u$,
then the following two statements are equivalent:
\begin{itemize}
\item $t=1$.
\item For all partitions $\lambda$ with $S_n(\lambda)\subseteq D$ and all permutations $\mu$, $\nu\in S_n(\lambda)$,
the following equality holds: 
$\frac{\textnormal{Pr}(\mu)}{\textnormal{Pr}(\nu)}=\frac{\textnormal{Pd}(\mu)}{\textnormal{Pd}(\nu)}$. That is, the
ratio between steady state probabilities does not change in moving from the {\ASEP} to the {\DASEP}.
\end{itemize}
\end{proposition}

\noindent In fact, we conjecture the following more general statement.
\begin{conjecture} \label{conj:name1}
If $D=\textnormal{{\DASEP}}(n,p,q)$ is parameterized as described above by $t$ and $u$,
then the following two statements are equivalent:
\begin{itemize}
\item $t=1$.
\item For all partitions $\lambda$ with $S_n(\lambda)\subseteq D$ and all permutations $\mu$, $\nu\in S_n(\lambda)$,
the following equality holds: 
$\frac{\textnormal{Pr}(\mu)}{\textnormal{Pr}(\nu)}=\frac{\textnormal{Pd}(\mu)}{\textnormal{Pd}(\nu)}$. That is, the
ratio between steady state probabilities does not change in moving from the {\ASEP} to the {\DASEP}.
\end{itemize}
\end{conjecture} 
\begin{proof}[Partial proof] We will prove this only in the $\implies$ direction. If $t=1$ we can replace $\lambda$ with
a similar partition but with species of the same type being replaced by similar distinct species. For example,
if $\lambda=(3,3,3,2,1,0,...)$ we would map this to $\hat{\lambda}=(3_1,3_2,3_3,2,1,0,...)$ and allow adjacent species
originally of the same type to be exchanged with the same transition probability. This will create a completely
symmetric situation, so all steady state probabilities are equal. As an equal number of $\hat{\lambda}$'s can be
derived from each $\lambda$ this means all original steady state probabilities are equal as well, so
$\frac{\textnormal{Pr}(\mu)}{\textnormal{Pr}(\nu)}=\frac{\textnormal{Pd}(\mu)}{\textnormal{Pd}(\nu)}=1$. This
completes the proof in the $\implies$ direction.
\end{proof}

\pagebreak
Let us motivate Conjecture~\ref{conj:name1}  by showing that it holds on one example.
Following are the nine values (of the nine steady state probabilities) we 
must solve for to prove this conjecture for {\DASEP}$(3,3,2)$:
\begin{align}  &a_1 = \textnormal{Pd}(0,1,1), \qquad  a_2 = \textnormal{Pd}(0,2,2), \qquad a_3 = \textnormal{Pd}(0,3,3),\\
 &b_1 = \textnormal{Pd}(0,2,3),  \qquad  b_2 = \textnormal{Pd}(0,1,3), \qquad b_3 = \textnormal{Pd}(0,1,2),\\
& c_1 = \textnormal{Pd}(0,3,2),  \qquad c_2 = \textnormal{Pd}(0,3,1),   \qquad c_3 = \textnormal{Pd}(0,2,1).\end{align}

These values can be obtained by solving the following set of nine equations:
\begin{equation}\begin{cases}
2ua_1=b_3+c_3 \\
 (2+t)b_3+ub_3+ub_3+b_3=(1+2t)c_3+b_2+a_2+ua_1 \\
(1+2t)c_3+uc_3+uc_3+c_3=(2+t)b_3+c_2+a_2+ua_1 \\
a_2+a_2+ua_2+ua_2=b_1+c_1+ub_3+uc_3 \\
(2+t)b_2+ub_2+b_2=(1+2t)c_2+b_1+ub_3 \\
 (1+2t)c_2+uc_2+c_2=(2+t)b_2+c_1+uc_3 \\
 2a_3=ub_1+uc_1 \\
 (2+t)b_1+ub_1+b_1+b_1=(1+2t)c_1+a_3+ub_2+ua_2 \\
 (1+2t)c_1+uc_1+c_1+c_1=(2+t)b_1+a_3+ua_2+uc_2. 
\end{cases}\end{equation}
Without working through all the details, this can be solved to give
\begin{equation}
\begin{split}
& (4u^3+36u^2t+90ut^2+72t^3+32u^2+206ut+270t^2+108u+322t+120)c_3 \\
& =(4u^3+24u^2t+54ut^2+36t^3+44u^2+190ut+198t^2+160u+350t+200)b_3.
\end{split}
\end{equation}

As previously discussed,
$\textnormal{Pr}(0,1,2) = \frac{1+2t}{9(1+t)}$ and
$\textnormal{Pr}(2,1,0) = \frac{2+t}{9(1+t)}$, so for $b_3=\textnormal{Pd}(0,1,2)$ and $c_3=\textnormal{Pd}(2,1,0)$
to be in the same ratio we would require $b_3=k(1+2t)$ and $c_3=k(2+t)$ for some $k$. It follows, 
after also dividing through by $2$, that
\begin{equation}
\begin{split}
& (2u^3+18u^2t+45ut^2+36t^3+16u^2+103ut  +135t^2+54u+161t+60)(t+2) \\
& =(2u^3+12u^2t+27ut^2+18t^3+22u^2+95ut  +99t^2+80u+175t+100)(2t+1).
\end{split}
\end{equation}
This can be expanded to
\begin{equation}
\begin{split}
& 2u^3t+18u^2t^2+45ut^3+36t^4+4u^3+52u^2t+193ut^2  \\
&+207t^3+32u^2+260ut+431t^2+108u+382t+120 \\
=\ & 4u^3t+24u^2t^2+54ut^3+36t^4+2u^3+56u^2t+\\
& 217ut^2  +216t^3+22u^2+255ut+449t^2+80u+375t+100.
\end{split}
\end{equation}
This can be reduced to
\begin{equation}
 2u^3t+6u^2t^2+9ut^3-2u^3+4u^2t+24ut^2  +9t^3-10u^2-5ut+18t^2-28u-7t-20=0.
\end{equation}
This can be factored as
\begin{equation}
\begin{split}
& (t-1)(2u^3+6u^2t+9ut^2+10u^2+33ut+9t^2+28u+27t+20)=0.
\end{split}
\end{equation}

Since $u\ge 0$ and $t\ge 0$, it follows that $t=1$. This completes the proof in the $\impliedby$ direction for the
{\DASEP}$(3,3,2)$ case.

\section{Proof of the conjecture for \texorpdfstring{{\DASEP}$(3,p,2)$}{DASEP(3,p,2)}}
It would be an endless game to prove the conjecture ``case by case'', with more and more cumbersome computations,
so let us now prove it for an infinite family of models.
More precisely, we now prove Conjecture~\ref{conj:name1} for {\DASEP}$(3,p,2)$
(our previous examples covered the cases $p=2$  and $p=3$). 
To solve this case we essentially need to solve for each of $p^2$ prior probabilities $p_{i,j}=\textnormal{Pd}(0,i,j)$
for $1\le i,j\le p$. The steady state probabilities can be obtained by solving a set of $p^2$ linear equations each of which
essentially demands equilibrium for each of the possible states of the process. The generic form of such an equation, for $i<j$, is
given by
\begin{equation}
\label{eq44}
(4+t+2u)p_{i,j}=(1+2t)p_{j,i}+p_{i+1,j}+p_{i,j+1}+up_{i-1,j}+up_{i,j-1}.
\end{equation}
For $i>j$ the equation is
\begin{equation}  (3+2t+2u)p_{i,j}=(2+t)p_{j,i}+p_{i+1,j}+p_{i,j+1}+up_{i-1,j}+up_{i,j-1} \end{equation}
For $i=j$ the equation simplifies to
\begin{equation}  (2+2u)p_{i,i}=p_{i+1,i}+p_{i,i+1}+up_{i-1,i}+up_{i,i-1}. \end{equation}

The equation may be similarly simplified for other edge cases such as $i=1<j$, $i<j=p$, $i=1<j=p$, $i>j=1$,
$i=p>j$, $i=p>j=1$, $i=j=1$, and $i=j=p$. For the sake of brevity we do not list all such cases in detail.

From the first above equation we can define a polynomial $A_{i,j}$ by gathering all terms on the left:
\begin{equation}  A_{i,j}:=(4+t+2u)p_{i,j}-(1+2t)p_{j,i}-p_{i+1,j}-p_{i,j+1}-up_{i-1,j}-up_{i,j-1}. \end{equation}

We can similarly define $A_{i,j}$ under the conditions stated for the various edge cases. We next define a $p^2\times p^2$ 
matrix $B$ as follows:
\begin{equation} B_{p(i_1-1)+j_1,p(i_2-1)+j_2}=[p_{i_1,j_1}] A_{i_2,j_2}. \end{equation}

The next step is to prove that the rank of $B$ is $p^2-1$. To see this, we first observe that the sum of all rows of $B$ is identically
zero, meaning that the rank cannot be $p^2$. For the rank to then be $p^2-1$, we would then need to show that no nontrivial linear 
combination of a proper subset of the rows can be zero. If we let row $i,j$ be $R_{i,j}$ and for some coefficients $c_{i,j}$ we have
$\sum_{i,j} c_{i,j}R_{i,j}=0$, then we need to show that if any $c_{i,j}=0$, then all $c_{i,j}=0$. The only rows with a $t$ term in
column $i,j$ will be $R_{i,j}$ and $R_{j,i}$. Hence if $c_{i,j}=0$, it follows that $c_{j,i}$=0.

We next show that if $c_{i,j}=0$ it follows that $c_{i-1,j-1}=0$. We can do this by first showing that $c_{i-1,j}$ and $c_{i,j-1}$
must be negations of one another. The only rows with a $u$ term in column $i,j$ will be $R_{i,j}$, $R_{i-1,j}$, and $R_{i,j-1}$, with 
the latter two having the same coefficient. Hence the following two statements are equivalent: $c_{i,j}=0$ and
$c_{i-1,j}+c_{i,j-1}=0$. We can similarly show that $c_{i,j}=0$ and $c_{i+1,j}+c_{i,j+1}=0$ are equivalent. So from $c_{i,j}=0$
we can derive $c_{i-1,j-1}=0$. By repeated application of the same argument we will get $c_{k,1}=0$ or $c_{1,k}=0$ for some $k$.
\pagebreak

Likewise, using the equations for the edge cases $i=1<j$ and $i>j=1$, the only rows with a $u$ term in column $1,k$ will be $R_{1,k}$ and $R_{1,k-1}$ and the only 
rows with a $u$ term in column $k,1$ will be $R_{k,1}$ and $R_{k-1,1}$. So from $c_{k,1}=0$ we can derive $c_{k-1,1}=0$ and from
$c_{1,k}=0$ we can derive $c_{1,k-1}=0$. By repeated application of this we will get to $c_{1,1}=0$. 
By reversing the above arguments
it follows that $c_{i,j}=0$ for any $i,j$ and we have proven:

\begin{lemma} \label{lemma}
The rank of the matrix $B$ as defined above is $p^2-1$. 
\end{lemma} 

We next prove a result about the values of the $p_{i,j}$.
\begin{proposition}One has 
\begin{equation}  p_{i,j} + p_{j,i} = \frac{2u^{i+j-2}}{\left(\sum_{k=0}^{n-1} u^k\right)^2} \text{\qquad and \qquad}  p_{i,i} = \frac{u^{2i-2}}{\left(\sum_{k=0}^{n-1} u^k\right)^2}. \end{equation}
\end{proposition}
\begin{proof} This can be proven by eliminating the variable $t$ from the set of linear equations above. For example, if
we add the equations for $i<j$ and $j<i$ we get the following:
\begin{equation}
\begin{split}
& (4+t+2u)p_{i,j}+(3+2t+2u)p_{j,i} \\
=\ & (2+t)p_{i,j}+(1+2t)p_{j,i}+p_{i+1,j}+p_{j,i+1}+p_{i,j+1}+p_{j+1,i} \\
& +up_{i-1,j}+up_{j,i-1}+up_{i,j-1}+up_{j-1,i}.
\end{split}
\end{equation}

If we let $q_{i,j}=p_{i,j}+p_{j,i}$ the above can be simplified to
\begin{equation}
\begin{split}
& (2+2u)q_{i,j}
=q_{i+1,j}+q_{i,j+1}
+uq_{i-1,j}+uq_{i,j-1}.
\end{split}
\end{equation}

If we substitute in the values for $q_{i,j}$ from the theorem we are attempting to prove to the above equation, we see that
it does satisfy the above equation. Therefore the values of $q_{i,j}$ given in the theorem represent one possible feasible solution
to the set of equations. Moreover, via Lemma~\ref{lemma} about the rank of $B$, the solution must be unique. This completes
the proof.\end{proof}

To continue with the proof of Conjecture \ref{conj:name1} in the $\impliedby$ direction, we note that from
$\frac{\textnormal{Pr}(\mu)}{\textnormal{Pr}(\nu)}=\frac{\textnormal{Pd}(\mu)}{\textnormal{Pd}(\nu)}$ it follows that
$\frac{\textnormal{Pr}(0,2,1)}{\textnormal{Pr}(0,1,2)}=\frac{\textnormal{Pd}(0,2,1)}{\textnormal{Pd}(0,1,2)}$ or
$\frac{2+t}{1+2t}=\frac{p_{2,1}}{p_{1,2}}$. This expands as $(2+t)p_{1,2}=(1+2t)p_{2,1}$. From the above theorem we know that
\begin{equation} p_{1,2} + p_{2,1} = \frac{2u}{\left(\sum_{k=0}^{n-1} u^k\right)^2}. \end{equation}
We can then solve for $p_{1,2}$ giving
\begin{equation}  p_{1,2} = \frac{2(1+2t)u}{3(1+t)\left(\sum_{k=0}^{n-1} u^k\right)^2}. \end{equation}
From the equation~\eqref{eq44} for $i=1<j$ we get
\begin{equation} (3+t+2u)p_{1,2}=(1+2t)p_{2,1}+p_{2,2}+p_{1,3}+up_{1,1}. \end{equation}
Substitute in to get
\begin{equation}
 \frac{2(3+t+2u)(1+2t)u}{3(1+t)\left(\sum_{k=0}^{n-1} u^k\right)^2 }
 = \frac{2(2+t)(1+2t)u}{3(1+t) \left(\sum_{k=0}^{n-1} u^k\right)^2 }
+\frac{3(1+t)(1+u)u}{3(1+t) \left(\sum_{k=0}^{n-1} u^k\right)^2}
+p_{1,3}.
\end{equation}  
This simplifies to
\begin{equation}
 p_{1,3} = \frac{(5ut+u+t-1)u}{3(1+t)\left(\sum_{k=0}^{n-1} u^k\right)^2}.
\end{equation}
A similar argument to that used to produce the above equation for $p_{1,2}$ will give us
\begin{equation} p_{1,3} = \frac{2(1+2t)u^2}{3(1+t)\left(\sum_{k=0}^{n-1} u^k\right)^2}. \end{equation}
Equating the last two equations and solving gives us $t=1$. This completes the proof of
\begin{theorem}
Conjecture \ref{conj:name1} holds for $D=\textnormal{{\DASEP}}(3,p,2)$.
\end{theorem}

\section{Future work}

Three main potential directions for future work are indicated. One is that further results should be obtained 
with a view to eventually proving Conjecture \ref{conj:name1}. We proved it for $\textnormal{{\DASEP}}(3,p,2)$
and the suggestion would be to prove it for $\textnormal{{\DASEP}}(n,2,2)$ and $\textnormal{{\DASEP}}(n,2,q)$ before
eventually proceeding to $\textnormal{{\DASEP}}(n,p,q)$. Similarly considering the case where $0$ represents a ball
with species $0$ rather than the absence of a species is a variant that should be explored.
The other, and more ambitious,
possible goal for future research would be to come up with a complete combinatorial characterization of the steady
state probabilities for the {\DASEP}. For the {\ASEP}, this has been done in~\cite{multiline} and~\cite{martin} leading to
a deep relationship being discovered between the {\ASEP} and Macdonald polynomials. 

\medskip

\subsection*{Acknowledgments}
Thanks to the editors of this special issue as well as to the anonymous referees for their many helpful
suggestions in the review and publication of this paper.
\smallskip
\subsection*{Funding} 
Thanks to Real Time Agents Inc of Pleasant Hill, CA, USA for partially funding this work.

\vspace{-1.2mm}

\bibliographystyle{SLC} 
\bibliography{ash}

\begin{thebibliography}{10}

\bibitem{ayyer}
Arvind Ayyer and Svante Linusson.
\newblock \href{https://arxiv.org/abs/1206.0316}{{An inhomogeneous multispecies
  TASEP on a ring}}.
\newblock  \textit{Advances in Applied Mathematics}, 57:21--43, 2014.

\bibitem{cantiniarc}
Luigi Cantini.
\newblock \href{https://arxiv.org/abs/1602.07921}{{Inhomogenous multispecies
  TASEP on a ring with spectral parameters}}.
\newblock  \textit{arXiv}, 2016.

\bibitem{cantini}
Luigi Cantini, Jan de~Gier, and Michael Wheeler.
\newblock \href{https://arxiv.org/abs/1505.00287}{{Matrix product formula for
  Macdonald polynomials}}.
\newblock  \textit{Journal of Physics A: Mathematical and Theoretical},
  48:38--64, 2015.

\bibitem{multiline}
Sylvie Corteel, Olya Mandelshtam, and Lauren Williams.
\newblock
  \href{http://math.univ-lyon1.fr/homes-www/slc/wpapers/FPSAC2019/97.pdf}{{From
  multiline queues to Macdonald polynomials via the exclusion process}}.
\newblock  \textit{Proceedings of the 31st Conference on Formal Power Series
  and Algebraic Combinatorics}, 82B:1--12, 2019.

\bibitem{ferrari}
Pablo Ferrari and James Martin.
\newblock \href{http://dx.doi.org/10.1214/009117906000000944}{Stationary
  distributions of multi-type totally asymmetric exclusion processes}.
\newblock  \textit{Annals of Probability}, 35:807--832, 2007.

\bibitem{kim}
Donghyun Kim and Lauren Williams.
\newblock \href{https://arxiv.org/abs/2102.00560}{{Schubert polynomials and the
  inhomogeneous TASEP on a ring}}.
\newblock  \textit{Séminaire Lotharingien de Combinatoire}, 20:1--12, 2021.

\bibitem{lam}
Thomas Lam and Lauren Williams.
\newblock \href{https://arxiv.org/abs/1102.4406}{{A Markov chain on the
  symmetric group that is Schubert positive?}}
\newblock  \textit{Experimental Mathematics}, 21:189--192, 2012.

\bibitem{macdonald}
Ian~G. Macdonald.
\newblock \href{https://math.berkeley.edu/~corteel/MATH249/macdonald.pdf}{
  \textit{{Symmetric Functions and Hall Polynomials}}}.
\newblock Oxford University Press, 1995.

\bibitem{mandelshtam}
Olya Mandelshtam.
\newblock \href{https://arxiv.org/abs/1707.02663}{{Toric tableaux and the
  inhomogeneous two-species TASEP on a ring}}.
\newblock  \textit{Advances in Applied Mathematics}, 113:101958, 2020.

\bibitem{martin}
James Martin.
\newblock \href{https://doi.org/10.1214/20-EJP421}{{Stationary distributions of
  the multi-type ASEP}}.
\newblock  \textit{Electronic Journal of Probability}, 43:1--41, 2020.

\end{thebibliography}
\end{document}